# Sensitivity integrals and related inequalities for process control systems

Shaival Nagarsheth*
Electrical Engineering Department, Sardar Vallabhbhai National Institute of Technology (SVNIT), Surat, Gujarat, India -395007. E-mail: shn411@gmail.com

Shambhu N Sharma
Electrical Engineering Department, Sardar Vallabhbhai National Institute of Technology (SVNIT), Surat, Gujarat, India -395007. E-mail: snsvolterra@gmail.com

**Abstract**

This paper exhibits the closed-loop design constraints using the non-analytic function theory. First, the paper generalizes the sensitivity integral for linear feedback systems with the non-analytic sensitivity function. Sensitivity inequalities are determined by the integral relationships based on the presence of non-minimum phase zeros and right half plane poles. These inequalities are rephrased in plant parameter context, which must be satisfied by the feedback design. That indicates the ability of controllers under the influence of input disturbances and plant parameter variations. The paper then extends the integral to the analytic sensitivity function of the augmented linear feedback systems. This is useful to augment the ability of a linear feedback system to handle input disturbances and plant uncertainties, via modified sensitivity function theory. Numerical simulations are carried out to perform sensitivity analysis on three chemical control systems. That describes the usefulness and demonstrates the applicability of the result of this paper to examine and augment the ability of a linear feedback system.

## 1. Introduction

The role of sensitivity functions is attributed to the ubiquity of control systems in the diverse field. The sensitivity function has found applications to study the qualitative characteristics of linear feedback systems subject to variation in parameters, stability margin as well as loop robustness (Xie and Lei, 2000). The goal in the design of the controller is not only to achieve the desired output but also loop robustness. Zames and Francis (1985) have explained the usefulness of sensitivity function; see Doyle and Stein (1981) as well. The importance of sensitivity function is directly linked to disturbance rejection, small parameter variations and tracking errors (Boyd and Desoer, 1985).

Rolf Nevanlinna introduced the Poisson-Jensen formula, which is a direct consequence of Jensen's formula (Jensen, 1899) by applying the Möbius transformation in the discrete domain. Subsequently, Levinson and Redheffer (1970) proved the same result in continuous settings as a particular case (Fisher, 1971). Freudenberg and Looze (1983) represented the Poisson-Jensen integral as a constraint on the feedback design of control systems. H W Bode (1945) studied the constraints on the linear time-invariant feedback system with single input single output settings. However, Horowitz (1963) interpreted Bode's work to control systems. The design tradeoffs and limitations of closed-loop feedback systems can be expressed in terms of the sensitivity function. Bode's Theorem states that for the stable open-loop transfer function, the log magnitude of the sensitivity function must equal to zero. This result of Bode's Theorem changes and the integral constraints vary with the presence of open-loop Non-Minimum Phase (NMP) zeros and Right Half Plane (RHP) poles (Freudenberg and Looze, 1985). A similar concept has been extended to linear multivariable feedback systems (Chen, 1995). The Bode integral was extended to linear time-varying systems (Iglesias, 2001) as well as linear time-periodic systems (Sandberg and Bernhardsson, 2005). This integral relation states that if the disturbance attenuation is improved in the sense of the sensitivity reduction in a frequency range, it worsens in another frequency range. This property of the constraint is described as a waterbed effect (Åström and Murray, 2008).

Researchers have attempted to apply sensitivity integrals to practical problems (Åström, 2000; Stein, 2003; Costa-Castelló and Dormido, 2015). Stein (2003) unfolded sensitivity integrals for aerodynamic problems. Sensitivity aspects to address the plant uncertainty in the ship steering system is with Roberts et al. (1988). Chu et al. (2018) have presented a sensitivity based graphical PID method to deal with uncertain industrial plants.

Chemical control systems are also subjected to input disturbances and plant parameter variations stemming from various sources (Åström, 2000). Control tuning achieving the robustness to address the uncertainty, makes a control problem hard (Bernstein, 2002). Such problems can be resolved by examining the controller's ability to handle uncertainty as well as augmenting the ability using formal methods. Examining as well as augmenting can be achieved using the sensitivity integral and the logarithmic sensitivity function. However, publications on filling the niche between sensitivity analysis and chemical control processes are sparingly cited (Vilanova and Alfaro, 2013; Tofighi et al., 2015). Although various controllers are proposed to achieve the robustness level, the analysis from sensitivity integral perspective to examine their controller ability needs to be explored further.

This paper attempts two control problems: (i) the problem of examining the ability of linear feedback systems to handle the effect of input disturbances and plant parameter variations (ii) the problem of augmenting the ability of linear feedback systems to handle the effect of input disturbances and plant parameter variations. The problem of examining the ability of linear feedback systems utilizes the non-analytic function theory. On the other hand, the problem of augmenting the ability of linear feedback systems adopts the analytic function theory. The ability is augmented by introducing a compensator in the existing forward path of the closed-loop system through a modified sensitivity function. In this paper, sensitivity analysis has been carried out revealing the controllers' ability to handle plant uncertainty in lieu of reference tracking. The controller might be suitable at reference tracking and give a stable response. But under plant uncertainty the controller might fail and may lead to an increase in the sensitivity of the closed-loop system, leading to instability. This paper focuses on 'integrating analytic function, complex integration, and sensitivity function.' Theorems 1 to 4 are constructed to attempt the above mentioned two control problems, with proofs, and are then applied to three chemical control systems. The sensitivity integrals, as well as their inequalities, become the design constraints for feedback system (Costa-Castelló and Dormido, 2015). These are useful in examining the ability of the feedback system under plant uncertainty and input disturbances. The usefulness of Sensitivity inequalities generated, using Poisson-Jenson and Bode's sensitivity integrals, is demonstrated by numerical simulations carried out in this paper.

2. **Main results**

In order to achieve the objective of the paper, first, the non-analytic function theory is used to examine the ability of linear feedback systems to handle the effect of input disturbances and plant parameter variations. That is formalized in Theorems 1 and 2. Then, they are applied to three chemical control systems. Furthermore, the analytic function theory is used to augment the ability to handle the effect of input disturbances and plant parameter variations. This can be found in Theorems 3 and 4 of the paper. Finally, Theorems 3 and 4 are applied to a chemical control example.

**2.1 Sensitivity integrals and related inequalities**

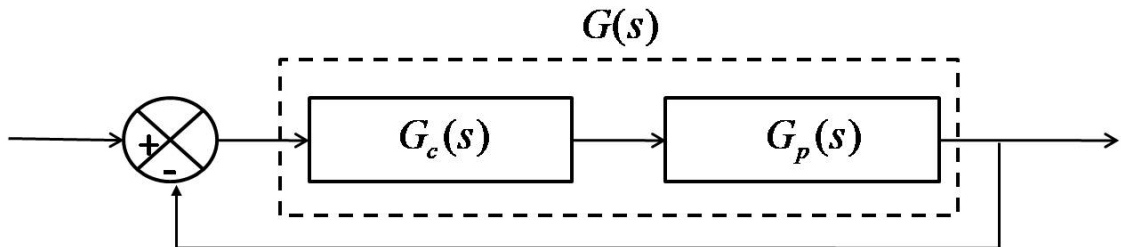

**Figure 1.** Linear feedback system

Figure 1 shows the block diagram of a linear feedback system taken into consideration. Here $G_p(s)$ is a plant transfer function, $G_c(s)$ is the controller transfer function and $G(s)$ is the forward path transfer function of the feedback loop. The sensitivity function $g(s)$ is defined as

$$g(s) = \frac{1}{1+G(s)}, \quad (1)$$

where $G(s) = G_c(s)G_p(s)$. In the structure of the sensitivity function $g(s)$ the NMP zeros of $g(s)$ can be regarded as unstable open-loop poles. The poles of the function $g(s)$ become the closed-loop poles of the linear feedback system. To cover a wider class of cases, we consider the non-analytic sensitivity function consisting of non-minimum phase zeros and right half plane poles. The sensitivity integral of the sensitivity function gives a trajectory of the absolute sensitivity over a range of frequencies for the linear feedback system.

**Theorem 1.** Consider a feedback control loop with the sensitivity function $g(s)$, which has the following properties:

(i) the function $g(s)$ is nonanalytic

(ii) the sensitivity function $g(s)$ has non-minimum phase zeros $\alpha_1, \alpha_2, \ldots, \alpha_M$ in the open Right Half Plane (RHP)

(iii) the sensitivity function $g(s)$ has $\beta_1, \beta_2, \ldots, \beta_N$ RHP poles. Then, the Poisson-Jensen formula becomes

$$\int_{-\infty}^{\infty} \log|g(j\omega)| \frac{\sigma_k}{\sigma_k^2 + (\omega - \eta_k)^2} d\omega = -\pi(-\log|g(\sigma_k + j\eta_k)| + \log \sum_{1 \leq i \leq M} \log\left|\frac{\sigma_k + j\eta_k - \alpha_i}{\sigma_k + j\eta_k + \overline{\alpha}_i}\right|$$

$$+ \sum_{1 \leq i \leq N} \log\left|\frac{\sigma_k + j\eta_k + \overline{\beta}_i}{\sigma_k + j\eta_k - \beta_i}\right|),$$

where $\sigma_k + j\eta_k$ is an arbitrary singular point in the RHP.

Suppose the open-loop transfer function $G(s)$ associated with the sensitivity function $g(s)$ has the NMP zeros. Consider the singular point $\sigma_k + j\eta_k = \varsigma_k$, where $\varsigma_k$ is the $k$th NMP open-loop zero associated with the sensitivity function $g(s)$. Then, the Poisson-Jensen formula becomes

$$\int_{-\infty}^{\infty} \log|g(j\omega)| \frac{\sigma_k}{\sigma_k^2 + (\omega - \eta_k)^2} d\omega = -\pi(\sum_{1 \leq i \leq M} \log\left|\frac{\sigma_k + j\eta_k - \alpha_i}{\sigma_k + j\eta_k + \overline{\alpha}_i}\right| + \sum_{1 \leq i \leq N} \log\left|\frac{\sigma_k + j\eta_k + \overline{\beta}_i}{\sigma_k + j\eta_k - \beta_i}\right|), \text{ where } 1 \leq k \leq L \text{ and}$$

$L$ is restricted to the number of NMP open-loop zeros.

**Proof.** We associate a singular point $s_0$ with the sensitivity function $g(s)$. The singular point $s_0$ is arbitrary or coincides with one of the NMP zeros of the open-loop transfers function $G(s)$. After associating the singular point, i.e., lying in the RHP, with the logarithmic function $\log g(s)$, the non-analyticity arises. The NMP zeros of the sensitivity function contained in the open RHP does not introduce the non-analyticity, but the logarithmic sensitivity function has the non-analyticity. We construct a new function $\widetilde{g}(s)$ with the analyticity, i.e.

$$\widetilde{g}(s) = g(s) \prod_{1 \leq i \leq M} \left(\frac{s + \overline{\alpha}_i}{s - \alpha_i}\right) \prod_{1 \leq i \leq N} \left(\frac{s - \beta_i}{s + \overline{\beta}_i}\right), \quad \frac{\widetilde{g}(s)}{s - s_0} = g(s) \frac{\prod_{1 \leq i \leq M} \left(\frac{s + \overline{\alpha}_i}{s - \alpha_i}\right) \prod_{1 \leq i \leq N} \left(\frac{s - \beta_i}{s + \overline{\beta}_i}\right)}{s - s_0},$$

(2a)

has the analyticity and non-analyticity, respectively. From the Cauchy integration Theorem (Levinson and Redheffer, 1970), we have

$$\log \widetilde{g}(s_0) = -\frac{1}{2\pi j} \oint \frac{\log \widetilde{g}(s)}{s - s_0} ds + \frac{1}{2\pi j} \oint \frac{\log \widetilde{g}(s)}{s - s_1}. \tag{2b}$$

Note that the point $s_1$, contained in the left half-plane, is an image of the singular point $s_0$. The function $\frac{\log \widetilde{g}(s)}{s - s_1}$ has the analyticity and the function $\frac{\log \widetilde{g}(s)}{s - s_0}$ has non-analyticity in the right half-plane, see Lemma A.1 of Freudenberg and Looze (1985). Furthermore,

$$\log \widetilde{g}(s_0) = \log g(s_0) + \sum_{1 \leq i \leq M} \log \frac{s_0 + \overline{\alpha}_i}{s_0 - \alpha_i} + \sum_{1 \leq i \leq N} \log \frac{s_0 - \beta_i}{s_0 + \overline{\beta}_i} = -\frac{1}{2\pi j} \oint \frac{(s_0 - s_1) \log \widetilde{g}(s)}{(s - s_0)(s - s_1)} ds.$$

The above is a consequence of equations (1) and (2) that can be further recast as

$$\log g(s_0) + \sum_{1 \leq i \leq M} \log \frac{s_0 + \overline{\alpha}_i}{s_0 - \alpha_i} + \sum_{1 \leq i \leq N} \log \frac{s_0 - \beta_i}{s_0 + \overline{\beta}_i} = -\frac{1}{2\pi j} \oint \frac{(s_0 - s_1)}{(s - s_0)(s - s_1)} (\log g(s)$$

$$+ \sum_{1 \leq i \leq M} \log \frac{s + \overline{\alpha}_i}{s - \alpha_i} + \sum_{1 \leq i \leq N} \log \frac{s - \beta_i}{s + \overline{\beta}_i}).$$

Thus,

$$\log g(s_0) + \sum_{1 \leq i \leq M} \log \frac{s_0 + \overline{\alpha}_i}{s_0 - \alpha_i} + \sum_{1 \leq i \leq N} \log \frac{s_0 - \beta_i}{s_0 + \overline{\beta}_i} = -\frac{1}{2\pi j} \oint \frac{(s_0 - s_1) \log g(s)}{(s - s_0)(s - s_1)} ds$$

$$-\frac{1}{2\pi j}\sum_{1\leq i\leq M}\oint\frac{s_0-s_1}{(s-s_0)(s-s_1)}\log\frac{s+\overline{\alpha}_i}{s-\alpha_i}ds$$

$$-\frac{1}{2\pi j}\sum_{1\leq i\leq N}\oint\frac{s_0-s_1}{(s-s_0)(s-s_1)}\log\frac{s-\beta_i}{s+\overline{\beta}_i}ds.$$

After decomposing integrals into parts, we arrive at

$$\log g(s_0)+\sum_{1\leq i\leq M}\log\frac{s_0+\overline{\alpha}_i}{s_0-\alpha_i}+\sum_{1\leq i\leq N}\log\frac{s_0-\beta_i}{s_0+\overline{\beta}_i}=-\frac{1}{2\pi j}\oint\frac{(s_0-s_1)\log g(s)}{(s-s_0)(s-s_1)}ds$$

$$-\frac{1}{2\pi j}\sum_{1\leq i\leq M}\oint\frac{s_0-s_1}{(s-s_0)(s-s_1)}\log\frac{s+\overline{\alpha}_i}{s-\alpha_i}ds$$

$$-\frac{1}{2\pi j}\sum_{1\leq i\leq N}\oint\frac{s_0-s_1}{(s-s_0)(s-s_1)}\log\frac{s-\beta_i}{s+\overline{\beta}_i}ds.$$

Furthermore,

$$\log g(s_0)+\sum_{1\leq i\leq M}\log\frac{s_0+\overline{\alpha}_i}{s_0-\alpha_i}+\sum_{1\leq i\leq N}\log\frac{s_0-\beta_i}{s_0+\overline{\beta}_i}=-\frac{1}{2\pi j}\int_{-\infty}^{\infty}\frac{j(s_0-s_1)\log g(j\omega)}{(j\omega-s_0)(j\omega-s_1)}d\omega$$

$$-\frac{1}{2\pi j}\underset{R\to\infty}{Lt}\int_{\frac{\pi}{2}}^{-\frac{\pi}{2}}\frac{(s_0-s_1)\log g(\mathrm{Re}^{j\gamma})}{(\mathrm{Re}^{j\gamma}-s_0)(\mathrm{Re}^{j\lambda}-s_1)}\mathrm{Re}^{j\gamma}\,jd\gamma$$

$$-\frac{1}{2\pi j}\sum_{1\leq i\leq M}\oint\frac{s_0-s_1}{(s-s_0)(s-s_1)}\log\frac{s+\overline{\alpha}_i}{s-\alpha_i}ds$$

$$-\frac{1}{2\pi j}\sum_{1\leq i\leq N}\oint\frac{s_0-s_1}{(s-s_0)(s-s_1)}\log\frac{s-\beta_i}{s+\overline{\beta}_i}ds.$$

Note that an open RHP can be regarded as a semi-circle with the radius $R\to\infty$ and

$$\underset{R\to\infty}{Lt}\frac{\log g(\mathrm{Re}^{j\gamma})}{(\mathrm{Re}^{j\gamma}-s_0)(\mathrm{Re}^{j\gamma}-s_1)}\mathrm{Re}^{j\gamma}=\underset{R\to\infty}{Lt}\frac{-\log(1+G(\mathrm{Re}^{j\gamma}))}{(\mathrm{Re}^{j\gamma}-s_0)(\mathrm{Re}^{j\gamma}-s_1)}\mathrm{Re}^{j\gamma}=0$$

for every relative degree of the open-loop transfer function $G(s)$.

Thus, $\log g(s_0)+\sum_{1\leq i\leq M}\log\frac{s_0+\overline{\alpha}_i}{s_0-\alpha_i}+\sum_{1\leq i\leq N}\log\frac{s_0-\beta_i}{s_0+\overline{\beta}_i}=-\frac{1}{2\pi j}\int_{-\infty}^{\infty}\frac{j(s_0-s_1)\log g(j\omega)}{(j\omega-s_0)(j\omega-s_1)}d\omega$

$$-\frac{1}{2\pi j}\sum_{1\leq i\leq M}\oint\frac{s_0-s_1}{(s-s_0)(s-s_1)}\log\frac{s+\overline{\alpha}_i}{s-\alpha_i}ds$$

$$-\frac{1}{2\pi j}\sum_{1\leq i\leq N}\oint\frac{s_0-s_1}{(s-s_0)(s-s_1)}\log\frac{s-\beta_i}{s+\overline{\beta}_i}ds. \quad (3)$$

Furthermore,

$$-\frac{1}{2\pi j}\sum_{1\leq i\leq M}\oint\frac{s_0-s_1}{(s-s_0)(s-s_1)}\log\frac{s+\overline{\alpha}_i}{s-\alpha_i}ds=-\frac{1}{2\pi j}\sum_{1\leq i\leq M}\int_{-\infty}^{\infty}\frac{j(s_0-s_1)}{(j\omega-s_0)(j\omega-s_1)}\log\frac{j\omega+\overline{\alpha}_i}{j\omega-\alpha_i}d\omega$$

$$-\frac{1}{2\pi j}\sum_{1\leq i\leq M}\underset{R\to\infty}{Lt}\int_{\frac{\pi}{2}}^{-\frac{\pi}{2}}\frac{(s_0-s_1)}{(\mathrm{Re}^{j\gamma}-s_0)(\mathrm{Re}^{j\gamma}-s_1)}\log(\frac{\mathrm{Re}^{j\gamma}+\overline{\alpha}_i}{\mathrm{Re}^{j\gamma}-\alpha_i})\mathrm{Re}^{j\gamma}\,jd\gamma$$

$$=-\frac{1}{2\pi j}\sum_{1\leq i\leq M}\int_{-\infty}^{\infty}\frac{j(s_0-s_1)}{(j\omega-s_0)(j\omega-s_1)}\log\frac{j\omega+\overline{\alpha}_i}{j\omega-\alpha_i}d\omega. \quad (4)$$

Under the limiting case, the second term of the right-hand side vanishes. The term

$$-\frac{1}{2\pi j}\sum_{1\leq i\leq N}\oint\frac{s_0-s_1}{(s-s_0)(s-s_1)}\log\frac{s-\beta_i}{s+\overline{\beta}_i}ds = -\frac{1}{2\pi j}\sum_{1\leq i\leq N}\int_{-\infty}^{\infty}\frac{j(s_0-s_1)}{(j\omega-s_0)(j\omega-s_1)}\log\frac{j\omega-\beta_i}{j\omega+\overline{\beta}_i}d\omega$$

$$-\frac{1}{2\pi j}\sum_{1\leq i\leq N}\underset{R\to\infty}{Lt}\int_{\frac{\pi}{2}}^{-\frac{\pi}{2}}\frac{(s_0-s_1)}{(Re^{j\gamma}-s_0)(Re^{j\gamma}-s_1)}\log(\frac{Re^{j\gamma}-\beta_i}{Re^{j\gamma}+\overline{\beta}_i})Re^{j\gamma}\,jd\gamma$$

$$=-\frac{1}{2\pi j}\sum_{1\leq i\leq N}\int_{-\infty}^{\infty}\frac{j(s_0-s_1)}{(j\omega-s_0)(j\omega-s_1)}\log\frac{j\omega-\beta_i}{j\omega+\overline{\beta}_i}d\omega. \quad (5)$$

Under the limiting case, the second term of the right-hand side of equation (5) vanishes as well. After combining equations (3)-(5), we get

$$\log g(s_0) + \sum_{1\leq i\leq M}\log\frac{s_0+\overline{\alpha}_i}{s_0-\alpha_i} + \sum_{1\leq i\leq N}\log\frac{s_0-\beta_i}{s_0+\overline{\beta}_i} = \frac{1}{2\pi j}\int_{-\infty}^{\infty}\frac{j(s_0-s_1)\log g(j\omega)}{(j\omega-s_0)(j\omega-s_1)}d\omega$$

$$-\frac{1}{2\pi j}\sum_{1\leq i\leq M}\int_{-\infty}^{\infty}\frac{j(s_0-s_1)}{(j\omega-s_0)(j\omega-s_1)}\log\frac{j\omega+\overline{\alpha}_i}{j\omega-\alpha_i}d\omega$$

$$-\frac{1}{2\pi j}\sum_{1\leq i\leq N}\int_{-\infty}^{\infty}\frac{j(s_0-s_1)}{(j\omega-s_0)(j\omega-s_1)}\log\frac{j\omega-\beta_i}{j\omega+\overline{\beta}_i}d\omega.$$

The idea is to simplify the above without approximating the expression. Here, we consider $s_0 = \sigma_k + j\eta_k$, $s_1 = -\sigma_k + j\eta_k$ as well as account for the real parts of both sides. As a result of this, we are led to a more straightforward expression, i.e.

$$\text{Re}(\log g(\sigma_k + j\eta_k) + \sum_{1\leq i\leq M}\log\frac{\sigma_k+j\eta_k+\overline{\alpha}_i}{\sigma_k+j\eta_k-\alpha_i} + \sum_{1\leq i\leq N}\log\frac{\sigma_k+j\eta_k-\beta_i}{\sigma_k+j\eta_k+\overline{\beta}_i}) = \text{Re}(\frac{1}{\pi}\int_{-\infty}^{\infty}\frac{\sigma_k\log g(j\omega)}{\sigma_k^2+(\omega-\eta_k)^2}d\omega$$

$$+\frac{1}{\pi}\sum_{1\leq i\leq M}\int_{-\infty}^{\infty}\frac{\sigma_k}{\sigma_k^2+(\omega-\eta_k)^2}\log\frac{j\omega+\overline{\alpha}_i}{j\omega-\alpha_i}d\omega$$

$$+\frac{1}{\pi}\sum_{1\leq i\leq N}\int_{-\infty}^{\infty}\frac{\sigma_k}{\sigma_k^2+(\omega-\eta_k)^2}\log\frac{j\omega-\beta_i}{j\omega+\overline{\beta}_i}d\omega).$$

Thus,

$$\log|g(\sigma_k+j\eta_k)| + \sum_{1\leq i\leq M}\log\left|\frac{\sigma_k+j\eta_k+\overline{\alpha}_i}{\sigma_k+j\eta_k-\alpha_i}\right| + \sum_{1\leq i\leq N}\log\left|\frac{\sigma_k+j\eta_k-\beta_i}{\sigma_k+j\eta_k+\overline{\beta}_i}\right| = \frac{1}{\pi}\int_{-\infty}^{\infty}\frac{\sigma_k\log|g(j\omega)|}{\sigma_k^2+(\omega-\eta_k)^2}d\omega$$

$$+\frac{1}{\pi}\sum_{1\leq i\leq N}\int_{-\infty}^{\infty}\frac{\sigma_k}{\sigma_k^2+(\omega-\eta_k)^2}\log\left|\frac{j\omega+\overline{\alpha}_i}{j\omega-\alpha_i}\right|d\omega$$

$$+\frac{1}{\pi}\sum_{1\leq i\leq M}\int_{-\infty}^{\infty}\frac{\sigma_k}{\sigma_k^2+(\omega-\eta_k)^2}\log\left|\frac{j\omega-\beta_i}{j\omega+\overline{\beta}_i}\right|d\omega.$$

The functions within the second and third integrals are vanishing. Finally,

$$\int_{-\infty}^{\infty}\log|g(j\omega)|\frac{\sigma_k}{\sigma_k^2+(\omega-\eta_k)^2}d\omega = -\pi(-\log|g(\sigma_k+j\eta_k)| + \sum_{1\leq i\leq M}\log\left|\frac{\sigma_k+j\eta_k-\alpha_i}{\sigma_k+j\eta_k+\overline{\alpha}_i}\right|$$

$$+ \sum_{1\leq i\leq N}\log\left|\frac{\sigma_k+j\eta_k+\overline{\beta}_i}{\sigma_k+j\eta_k-\beta_i}\right|). \quad (6)$$

After embedding $s_0 = \sigma_k + j\eta_k = \varsigma_k$, $s_1 = -\sigma_k + j\eta_k = -\overline{\varsigma}_k$, then we have

$$\int_{-\infty}^{\infty}\log|g(j\omega)|\frac{\sigma_k}{\sigma_k^2+(\omega-\eta_k)^2}d\omega = -\pi\left(\sum_{1\le i\le M}\log\left|\frac{\sigma_k+j\eta_k-\alpha_i}{\sigma_k+j\eta_k+\overline{\alpha}_i}\right| + \sum_{1\le i\le N}\log\left|\frac{\sigma_k+j\eta_k+\overline{\beta}_i}{\sigma_k+j\eta_k-\beta_i}\right|\right). \tag{7a}$$

The subscript $k$ is associated with the arbitrary singular point. For the case, the singular point coincides with the NMP open-loop zero, then the number of sensitivity integrals assuming the structure of equation (7a) would be the number of NMP open-loop zeros.

*Remark 1:* Suppose the sensitivity function $g(s)$ has no RHP poles, i.e., no RHP poles associated with the linear feedback system. The resulting sensitivity integral equation becomes

$$\int_{-\infty}^{\infty}\log|g(j\omega)|\frac{\sigma_k}{\sigma_k^2+(\omega-\eta_k)^2}d\omega = -\pi\sum_{1\le i\le M}\log\left|\frac{\sigma_k+j\eta_k-\alpha_i}{\sigma_k+j\eta_k+\overline{\alpha}_i}\right|. \tag{7b}$$

The above logarithmic integral in equation (7b) is a special case of the logarithmic integral in equation (7a) of Theorem 1. The sensitivity integral of Theorem 1 can be rephrased using the Bode integral. This can be found in case 1 of the *Appendix* of the paper.

**Theorem 2.** Suppose the sensitivity function obeys the properties of Theorem 1 of the paper. Suppose $\rho = \max_{0\le\omega\le\omega_c}|g(j\omega)| \in [0,1)$, $s_{max} = \max_{\omega_c<\omega}|g(j\omega)| \in (1,\infty)$. Then, the sensitivity inequalities are

$$\log s_{max} \ge \frac{\left|\frac{\pi}{2}(-\log|g(\sigma_k)| + \sum_{1\le i\le M}\log\left|\frac{\sigma_k-\alpha_i}{\sigma_k+\overline{\alpha}_i}\right| + \sum_{1\le i\le N}\log\left|\frac{\sigma_k+\overline{\beta}_i}{\sigma_k-\beta_i}\right|) + \tan^{-1}\frac{\omega_c}{\sigma_k}\log\rho\right|}{\frac{\pi}{2} - \tan^{-1}\frac{\omega_c}{\sigma_k}},$$

where $s_0 = \sigma_k$, $\eta_k = 0$, $-\log|g(\sigma_k)| + \sum_{1\le i\le M}\log\left|\frac{\sigma_k-\alpha_i}{\sigma_k+\overline{\alpha}_i}\right| + \sum_{1\le i\le N}\log\left|\frac{\sigma_k+\overline{\beta}_i}{\sigma_k+\beta_i}\right| < 0$, and

$$\log s_{max} \ge \frac{\left|\frac{\pi}{2}(\sum_{1\le i\le M}\log\left|\frac{\sigma_k-\alpha_i}{\sigma_k+\overline{\alpha}_i}\right| + \sum_{1\le i\le N}\log\left|\frac{\sigma_k+\overline{\beta}_i}{\sigma_k-\beta_i}\right|) + \tan^{-1}\frac{\omega_c}{\sigma_k}\log\rho\right|}{\frac{\pi}{2} - \tan^{-1}\frac{\omega_c}{\sigma_k}},$$

where $s_0 = \sigma_k = \varsigma_k, \eta_k = 0$, $\sum_{1\le i\le M}\log\left|\frac{\sigma_k-\alpha_i}{\sigma_k+\overline{\alpha}_i}\right| + \sum_{1\le i\le N}\log\left|\frac{\sigma_k+\overline{\beta}_i}{\sigma_k+\beta_i}\right| < 0$.

**Proof.** Consider the Poisson-Jensen sensitivity integral

$$\int_{-\infty}^{\infty}\log|g(j\omega)|\frac{\sigma_k}{\sigma_k^2+(\omega-\eta_k)^2}d\omega = -\pi\left(-\log|g(\sigma_k+j\eta_k)| + \sum_{1\le i\le M}\log\left|\frac{\sigma_k+j\eta_k-\alpha_i}{\sigma_k+j\eta_k+\overline{\alpha}_i}\right|\right.$$
$$\left.+ \sum_{1\le i\le N}\log\left|\frac{\sigma_k+j\eta_k+\overline{\beta}_i}{\sigma_k+j\eta_k+\beta_i}\right|\right).$$

Note that the integrand of the left-hand side is not an even function. To achieve a simplified and appealing form of the sensitivity integral, consider the singular points are the real. Then, the integrand is an even function. As a result of this, we arrive at

$$\int_0^{\infty}\log|g(j\omega)|\frac{\sigma_k}{\sigma_k^2+\omega^2}d\omega = -\frac{\pi}{2}\left(-\log|g(\sigma_k)| + \sum_{1\le i\le M}\log\left|\frac{\sigma_k-\alpha_i}{\sigma_k+\overline{\alpha}_i}\right| + \sum_{1\le i\le N}\log\left|\frac{\sigma_k+\overline{\beta}_i}{\sigma_k+\beta_i}\right|\right).$$

Then

$$\tan^{-1}\frac{\omega_c}{\sigma_k}\log\rho + \left(\frac{\pi}{2}-\tan^{-1}\frac{\omega_c}{\sigma_k}\right)\log s_{max} \ge -\frac{\pi}{2}\left(-\log|g(\sigma_k)| + \sum_{1\le i\le M}\log\left|\frac{\sigma_k-\alpha_i}{\sigma_k+\overline{\alpha}_i}\right|\right.$$
$$\left.+ \sum_{1\le i\le N}\log\left|\frac{\sigma_k+\overline{\beta}_i}{\sigma_k-\beta_i}\right|\right). \tag{8}$$

After rearranging the terms, embedding the conditions, we arrive at two sensitivity inequalities of Theorem 2. The sensitivity equality holds for frequency-independent transfer characteristics. On the other hand, equality will not hold for frequency-dependent transfer characteristics. The sensitivity inequality of Theorem 2 reveals the fact that logarithmic inequality arises from linear feedback systems.

*Remark 2:* Furthermore, the sensitivity inequality

$$\log s_{\max} \geq \frac{\left|\frac{\pi}{2}\sum_{1\leq i\leq M}\log\left|\frac{\sigma_k - \alpha_i}{\sigma_k + \overline{\alpha}_i}\right| + \tan^{-1}\frac{\omega_c}{\sigma_k}\log\rho\right|}{\frac{\pi}{2} - \tan^{-1}\frac{\omega_c}{\sigma_k}}, \quad (9)$$

where $\sum_{1\leq i\leq M}\log\left|\frac{\sigma_k - \alpha_i}{\sigma_k + \overline{\alpha}_i}\right| < 0$, $\beta_i = 0, \sigma_k = \varsigma_k, \eta_k = 0$. The sensitivity inequality in bound (9) is restrictive and a special case of (8). That is useful to achieve the sensitivity bound for the stable linear feedback system, where the system has a real open-loop non-minimum phase zeros. Notably, the sensitivity inequality in (8) covers a broader class of linear feedback systems. The restriction on the location of the sensitivity function's poles and zeros decides the sensitivity inequality. It is natural to construct the sensitivity inequality in the Bode terminology as well. This can be found in case 2 of the *Appendix* of the paper.

## 2.2 Modified sensitivity integrals and related inequalities

The problem of improving loop robustness reduces to the problem of reframing sensitivity bounds, which eventually leads us to the loop shaping problem for unstable open-loop systems (Åström and Murray, 2008). The system's sensitivity gets affected due to controller failure in the case of plant uncertainty and input disturbances. By developing a modified sensitivity function theory, a compensator is added into the existing forward path of the linear feedback system. This augments the ability of linear feedback system to handle process parameter variations as well as input disturbances. The modified sensitivity integrals and related inequalities are summarized in Theorems 3 and 4 of the paper.

**Theorem 3.** Consider the sensitivity function $g(s)$ obeys the properties of Theorem 1. The modified sensitivity function $\widetilde{g}(s)$ is analytic, $\widetilde{g}(s) = g(s) \prod_{1\leq i\leq M}(\frac{s+\overline{\alpha}_i}{s-\alpha_i}) \prod_{1\leq i\leq N}(\frac{s-\beta_i}{s+\overline{\beta}_i})$. Then the Poisson-Jensen formula becomes

$$\int_{-\infty}^{\infty}\log|\widetilde{g}(j\omega)|\frac{\sigma_k}{\sigma_k^2 + (\omega - \eta_k)^2}d\omega = \pi\log|\widetilde{g}(\sigma_k + j\eta_k)|,$$

where $\sigma_k + j\eta_k$ is an arbitrary singular point in the RHP. Furthermore,

$$\int_{-\infty}^{\infty}\log|\widetilde{g}(j\omega)|\frac{\sigma_k}{\sigma_k^2 + (\omega - \eta_k)^2}d\omega = 0, 1 \leq k \leq L$$

for the open-loop transfer function $\widetilde{G}(s)$ associated with the sensitivity function $\widetilde{g}(s)$, where $L$ denotes NMP open-loop zeros. Suppose the singular point $\sigma_k + j\eta_k = \varsigma_k$, then

$$\int_{-\infty}^{\infty}\log|\widetilde{g}(j\omega)|\frac{\sigma_k}{\sigma_k^2 + (\omega - \eta_k)^2}d\omega = \pi\log|\widetilde{g}(\sigma_k + j\eta_k)|, 1 \leq k \leq S.$$

Note that $\varsigma_k$ is an NMP zero associated with the open-loop transfer $G(s)$ in lieu of the open-loop transfer function $\widetilde{G}(s)$. The integers $L$ and $S$ depend on the structure of $\widetilde{G}(s)$ and $G(s)$ respectively.

**Proof.** Consider the modified sensitivity function-sensitivity function relation

$$\widetilde{g}(s) = \kappa(s)g(s), \quad \kappa(s) = \prod_{1\leq i\leq M}(\frac{s+\overline{\alpha}_i}{s-\alpha_i}) \prod_{1\leq i\leq N}(\frac{s-\beta_i}{s+\overline{\beta}_i}).$$

Following equation (1a), we get

$$\frac{1}{1+\widetilde{G}(s)} = \frac{\prod_{1\leq i\leq M}(\frac{s+\overline{\alpha}_i}{s-\alpha_i}) \prod_{1\leq i\leq N}(\frac{s-\beta_i}{s+\overline{\beta}_i})}{1+G(s)}.$$

Then we have the following relation in the notational form,
$$\widetilde{G}(s) = \frac{1 - \kappa(s) + G(s)}{\kappa(s)}. \tag{10}$$

Utilizing equation (2a) of the paper,

$$\log \widetilde{g}(s_0) = -\frac{1}{2\pi j} \oint \frac{\log \widetilde{g}(s)}{s - s_0} ds + \frac{1}{2\pi j} \oint \frac{\log \widetilde{g}(s)}{s - s_1}$$

$$= -\frac{1}{2\pi j} \int_{-\infty}^{\infty} \frac{j(s_0 - s_1) \log \widetilde{g}(j\omega)}{(j\omega - s_0)(j\omega - s_1)} d\omega$$

$$-\frac{1}{2\pi j} \underset{R \to \infty}{Lt} \int_{\frac{\pi}{2}}^{-\frac{\pi}{2}} \frac{(s_0 - s_1) \log \widetilde{g}(Re^{j\gamma})}{(Re^{j\gamma} - s_0)(Re^{j\lambda} - s_1)} Re^{j\gamma} j d\gamma$$

$$= -\frac{1}{2\pi j} \int_{-\infty}^{\infty} \frac{j(s_0 - s_1) \log \widetilde{g}(j\omega)}{(j\omega - s_0)(j\omega - s_1)} d\omega.$$

Note that $\underset{R \to \infty}{Lt} \frac{\log \widetilde{g}(Re^{j\gamma})}{(Re^{j\gamma} - s_0)(Re^{j\gamma} - s_1)} Re^{j\gamma} = \underset{R \to \infty}{Lt} \frac{-\log(1 + \widetilde{G}(Re^{j\gamma}))}{(Re^{j\gamma} - s_0)(Re^{j\gamma} - s_1)} Re^{j\gamma} = 0$ for every relative degree of the open-loop transfer function $\widetilde{G}(s)$. Here, we consider $s_0 = \sigma_k + j\eta_k$, $s_1 = -\sigma_k + j\eta_k$ as well as account for the real parts of both sides. As a result of this, we are led to a simpler expression. Thus,

$$\log|\widetilde{g}(\sigma_k + j\eta_k)| = \frac{1}{\pi} \int_{-\infty}^{\infty} \frac{\sigma_k \log|\widetilde{g}(j\omega)|}{\sigma_k^2 + (\omega - \eta_k)^2} d\omega. \tag{11}$$

Consider the open-loop transfer function $\widetilde{G}(s)$ having NMP zeros, say the $k$th NMP zero of $\widetilde{G}(s)$ is $\varsigma_k$. Suppose the singular point $\sigma_k + j\eta_k = \varsigma_k$, then we get

$$\int_{-\infty}^{\infty} \log|\widetilde{g}(j\omega)| \frac{\sigma_k}{\sigma_k^2 + (\omega - \eta_k)^2} d\omega = 0, 1 \le k \le L. \tag{12}$$

The term $\log \widetilde{g}(j\omega)$ associated with Theorem 3 of the paper is analytic in the right half plane. On the other hand, the term $\log g(j\omega)$ associated with Theorem 1 is non-analytic.

Similar arguments can be made to arrive at the third sensitivity integral of Theorem 3 of the paper. Equation (10) in combination with the sensitivity integrals of Theorem 3, answers to the problem of augmenting the ability of linear feedback systems to handle input disturbances as well as plant parameter variations. Case 3 of the *Appendix* represents the bode terminology of the sensitivity integral of Theorem 3. 'Reshaping the forward path transfer function' using equation (10) and 'reshaping the sensitivity inequality' can be achieved using results of Theorem 3, which is sketched in Theorem 4 of the paper.

**Theorem 4.** Suppose the sensitivity function obeys the properties of Theorem 3 of the paper. Suppose $\rho = \underset{0 \le \omega \le \omega_c}{\max} |\widetilde{g}(j\omega)| \in [0,1)$, $s_{\max} = \underset{\omega_c < \omega}{\max} |\widetilde{g}(j\omega)| \in [1, \infty)$. Then, the sensitivity inequalities are

$$\log s_{\max} \ge \frac{\left| \frac{\pi}{2} \log|\widetilde{g}(\sigma_k)| - \tan^{-1} \frac{\omega_c}{\sigma_k} \log \rho \right|}{\frac{\pi}{2} - \tan^{-1} \frac{\omega_c}{\sigma_k}},$$

where $s_k = \sigma_k, \eta_k = 0, -\log|\widetilde{g}(\sigma_k)| < 0$. Furthermore,

$$\log s_{\max} \ge \frac{\left| \tan^{-1} \frac{\omega_c}{\sigma_k} \log \rho \right|}{\frac{\pi}{2} - \tan^{-1} \frac{\omega_c}{\sigma_k}},$$

where $s_k = \sigma_k = \varsigma_k, \eta_k = 0$ and $\varsigma_k$ is the $k$th real-valued NMP zero of the $\widetilde{G}(s)$.

**Proof.** Consider the Poisson-Jensen sensitivity integral for the sensitivity function $\widetilde{g}(s)$,

$$\int_{-\infty}^{\infty} \log|\widetilde{g}(j\omega)| \frac{\sigma_k}{\sigma_k^2 + (\omega - \eta_k)^2} d\omega = \pi \log|\widetilde{g}(\sigma_k + j\eta_k)|.$$

Consider the singular points are real, then the integrand is an even function. As a result of this, we arrive at

$$\int_{0}^{\infty} \log|\widetilde{g}(j\omega)| \frac{\sigma_k}{\sigma_k^2 + \omega^2} d\omega = \frac{\pi}{2} \log|\widetilde{g}(\sigma_k)|.$$

Since the closed-form expression of the logarithmic absolute sensitivity satisfying the above integral equation is not possible, we introduce the notion of the logarithmic inequality. As a result of this, we arrive at

$$\tan^{-1}\frac{\omega_c}{\sigma_k} \log \rho + (\frac{\pi}{2} - \tan^{-1}\frac{\omega_c}{\sigma_k}) \log s_{\max} \geq \frac{\pi}{2} \log|\widetilde{g}(\sigma_k)|.$$

After rearranging the terms, embedding the conditions, we arrive at two sensitivity inequalities. Construction of the sensitivity inequality in the Bode terminology for the modified sensitivity function $\widetilde{g}(s)$ can be found in case 4 of the *Appendix*.

## 3. Numerical simulations

From a design perspective, the use of a loop analysis tool is powerful. The impact of changes in the design of the controller $G_c(s)$ can be witnessed or measured if the bounds are in terms of the loop transfer function. The efficacy of a controller can be measured in several terms. The major one is the robustness to plant parameter variations and responses to reference signals and disturbances. The first step for a controller to provide adequate robustness performance is to satisfy the sensitivity bound. Practical examples of chemical control systems have frequency-dependent transfer characteristics. Hence, equality will not hold in the sensitivity bounds. The sensitivity function represents the disturbance attenuation and also relates the uncertainties in the plant parameter variations. Here, we take *three* examples with open-loop NMP zero and unstable poles. Generally, the usefulness of linear feedback systems is adjudged using integral square error, integral absolute error, and total variation. However, they do not suggest formally the efficacy of controllers in the sense of sensitivity. In this paper, we consider three chemical control systems and their controllers to examine their efficacy from the viewpoint of sensitivity by adopting the sensitivity inequality framework of Theorems of the paper.

**Example 1: FOIPDT (First Order Integrating Plus Dead Time)**

Consider a unity feedback system with the plant transfer function $G_p(s)$ (Bernstein and Hoagg, 2007), the controller transfer function $G_c(s)$ and the sensitivity function $g(s)$. The function $G_p(s)$ has NMP zeros and no RHP poles, i.e.

$$G_p(s) = \frac{k_p(1 - \tau_0 s)e^{-t_d s}}{s(1 + \tau s)}.$$

Consider the control problem of a benchmark system possessing a non-minimum phase behavior, e.g., a boiler-drum level system (Luyben, 2003). The boiler-drum level system is a subsystem that precedes the industrial reactor in chemical control systems. That is useful for heating chemical compounds (Åström and Bell, 2000).
To test the efficacy of controllers from the sensitivity-theoretic viewpoint for the system, consider the plant transfer function $G_p(s)$, the Luyben controller transfer function $G_{c1}(s)$ (Luyben, 2003) and the controller transfer function $G_{c2}(s)$ of Pai et al. (2010),

$$G_p(s) = \frac{0.547(1 - 0.418s)e^{-0.1s}}{s(1.06s + 1)}, \quad G_{c1}(s) = \frac{k_c(1 + \tau_I s)(1 + \tau_D s)}{\tau_I s(\alpha \tau_D s + 1)}, \quad G_{c2}(s) = k_c(1 + \frac{1}{\tau_I s} + \tau_D s). \tag{13}$$

Further, we adopt the sensitivity function notations $g_1(j\omega)$ and $g_2(j\omega)$ for the Luyben (2003) and Pai et al. (2010) controllers, respectively. The sensitivity function for both the control loops can be written as $g_i(s) = \frac{1}{1 + G_i(s)}$, $i = 1, 2$, where $1 \leq i \leq 2$, $G_i(s) = G_p(s)G_{ci}(s)$. Alternatively,

$$g_1(j\omega) = \frac{-\tau_I\omega^2(1-\alpha\tau_D\tau\omega^2) - j\tau_I(\alpha\tau_D+\tau)\omega^3}{-\tau_I\omega^2(1-\alpha\tau_D\tau\omega^2) - j(\alpha\tau_D+\tau)\tau_I\omega^3 + k_p k_c(a_1\cos t_d\omega - j\sin t_d\omega)(a_1+jb_1)}$$

$$= \frac{-x-jy}{k_p k_c(a_1\cos t_d\omega + b_1\sin t_d\omega) - x + j(k_p k_c(b_1\cos t_d\omega - a_1\sin t_d\omega) - y)}, \quad (14)$$

where $a_1 = (\tau_I\tau_0 + \tau_D\tau_0 - \tau_I\tau_D)\omega^2 + 1$, $b_1 = \tau_I\tau_D\tau_0\omega^3 + (\tau_I + \tau_D - \tau_0)\omega$, $x = \tau_I\omega^2(1-\alpha\tau_D\tau\omega^2)$, $y = \tau_I(\alpha\tau_D + \tau)\omega^3$. Moreover,

$$g_2(j\omega) = \frac{-\omega^2(\tau_I + j\omega\tau^2)}{-\omega^2(\tau_I + j\omega\tau^2) + k_p k_c((a_2\cos t_d\omega + b_2\sin t_d\omega) + j(b_2\cos t_d\omega - a_2\sin t_d\omega))}$$

$$= \frac{-\omega^2(\tau_I + j\omega\tau^2)}{k_p k_c(a_2\cos t_d\omega + b_2\sin t_d\omega) - \tau_I\omega^2 + j(k_p k_c(b_2\cos t_d\omega - a_2\sin t_d\omega) - \tau^2\omega^3)}, \quad (15)$$

where $a_2 = \tau_0(\tau_I - \tau_D)\omega^2 + 1$ and $b_2 = \tau_I\tau_D\tau_0\omega^3 + (\tau_I - \tau_0)\omega$.

The pole-zero property of the sensitivity function of the Luyben controller-embedded linear feedback system are the following: (i) the sensitivity function has no NMP zeros (ii) the sensitivity function has no RHP poles. Consider the singular point $s_0$, associated with the sensitivity integral, which coincides with the NMP open-loop zero, $s_0 = \frac{1}{\tau_0}$. As a result of this, we are led to the following sensitivity integral for the Luyben controller-embedded linear feedback system:

$$\int_{-\infty}^{\infty} \log|g_1(j\omega)| \frac{\tau_0}{1+\omega^2\tau_0^2} d\omega = 0 \quad (16)$$

with

$$|g_1(j\omega)| = \frac{\tau_I\omega^2\sqrt{(\alpha^2\tau_D^2\omega^2+1)(\tau^2\omega^2+1)}}{\sqrt{x^2+y^2+k_p^2 k_c^2(a_1^2+b_1^2) - 2k_p k_c((a_1 x + b_1 y)\cos t_d\omega + (b_1 x - a_1 y)\sin t_d\omega)}}.$$

The sensitivity integral for the Luyben controller-embedded linear feedback system is a consequence of Theorem 1 of the paper. The interpretation of the parameters $a_1, b_1, x$ and $y$ can be found in equation (14). Since the closed-form expression of the logarithmic absolute sensitivity associated with the sensitivity integral equation (16) is not possible, we achieve the Luyben controller sensitivity inequality. After adopting appropriate variables for the specific case $g_1(j\omega)$, Theorem 2 of the paper becomes the logarithmic inequality for the Luyben controller. Thus, we arrive at

$$\log s_{\max} > \frac{\left|\log \max_{0\leq\omega\leq\omega_c} |g_1(j\omega)| \tan^{-1}\omega_c\tau_0\right|}{\frac{\pi}{2} - \tan^{-1}\omega_c\tau_0} \quad (17)$$

in the Poisson-Jensen setting. The sensitivity inequality, for example 1 can also be written in the Bode inequality setting, i.e. $\log s_{\max} > \underset{\omega_l\to\infty}{Lt} \frac{|\omega_c \log\rho|}{\omega_l - \omega_c}$, see case 1 of the *Appendix* of the paper.

With similar pole-zero property, the logarithmic absolute sensitivity integral associated with the Pai et al. (2010) controller becomes

$$\int_{-\infty}^{\infty} \log|g_2(j\omega)| \frac{\tau_0}{1+\omega^2\tau_0^2} d\omega = 0 \quad (18)$$

with

$$|g_2(j\omega)| = \frac{\omega^2\sqrt{\omega^2\tau^4 + \tau_I^2}}{\sqrt{\omega^4(\omega^2\tau^4+\tau_I^2) + k_p^2 k_c^2(a_2^2+b_2^2) - 2k_p k_c\omega^2((b_2\omega\tau^2 + a_2\tau_I)\cos t_d\omega + (b_2\tau_I - a_2\omega\tau^2)\sin t_d\omega)}}.$$

For the notational brevity and convenience, we choose the above structure of the absolute sensitivity $|g_2(j\omega)|$ and the parameters $a_2$ and $b_2$ can be found in equation (15). After adopting appropriate variables for the specific case $g_2(j\omega)$, Theorem 2 of the paper becomes the logarithmic inequality for the Pai et al. (2010) controller. Thus, we get

$$\log s_{\max} > \frac{\left|\log \max_{0 \leq \omega \leq \omega_c} |g_2(j\omega)| \tan^{-1} \omega_c \tau_0\right|}{\frac{\pi}{2} - \tan^{-1} \omega_c \tau_0} \quad (19)$$

in the Poisson-Jensen setting.

Figure 2 shows a contrast between the graphical trajectories of the sensitivity functions $g_1(j\omega)$ and $g_2(j\omega)$ of the system $G_p(s)$ associated with the controllers $G_{c1}(s)$ and $G_{c2}(s)$ respectively. The controller parameters and sensitivity performance indices for both the control loops are listed in Table 1. The theory of the paper suggests that the sensitivity inequalities in (17) and (19) must be obeyed. In control literature, the Poisson-Jensen and Bode inequalities are significant logarithmic inequalities. The inequalities $\log s_{\max} > 9.48 \times 10^{-3}$ and $\log s_{\max} > 2.894 \times 10^{-5}$ hold for the Poisson-Jensen and Bode settings respectively for the sensitivity function $g_1(j\omega)$. For the sensitivity function $g_1(j\omega)$, the Poisson-Jensen's sensitivity inequality gives a refined lower bound than that of the Bode's inequality. For the sensitivity function $g_2(j\omega)$, the inequalities $\log s_{\max} > 0.0685$ and $\log s_{\max} > 2.077 \times 10^{-4}$ hold for the Poisson-Jensen and Bode settings, respectively. A calculation of the logarithmic absolute maximum sensitivity using the logarithmic sensitivity integral setup is intractable theoretically. We address the difficulty by adopting the graphical calculation and the logarithmic inequality interpretation. The sensitivity lower bounds for both the controllers are calculated analytically, and then we associate them with the logarithmic absolute maximum sensitivity. Finally, we get the logarithmic inequality. The graphical calculation of the logarithmic absolute maximum sensitivity obeys the logarithmic inequality, see Table 1. <span style="color:red">The numerical simulation depicted in Figure 2 concords the analytic calculation of the bounds in (17) and (19).</span> The bounds result from Theorem 2 of the paper.

Figure 2 reveals that the disturbances occurring at the frequencies such that $|g_i(j\omega)| < 1$ and $\log|g_i(j\omega)| < 0$ are attenuated, but the disturbances with frequencies such that $|g_i(j\omega)| > 1$ and $\log|g_i(j\omega)| > 0$ are amplified by the feedback. The graphical revelation in Figure 2 indicates that the Pai et al. (2010) controller has more attenuation than that of the Luyben (2003) controller for the disturbances in the frequency range $0 < \omega < \omega_c$. However, the disturbances with frequencies higher than $\omega_c$ are relatively less amplified in the controller $G_{c1}(s)$ in contrast to the controller $G_{c2}(s)$.

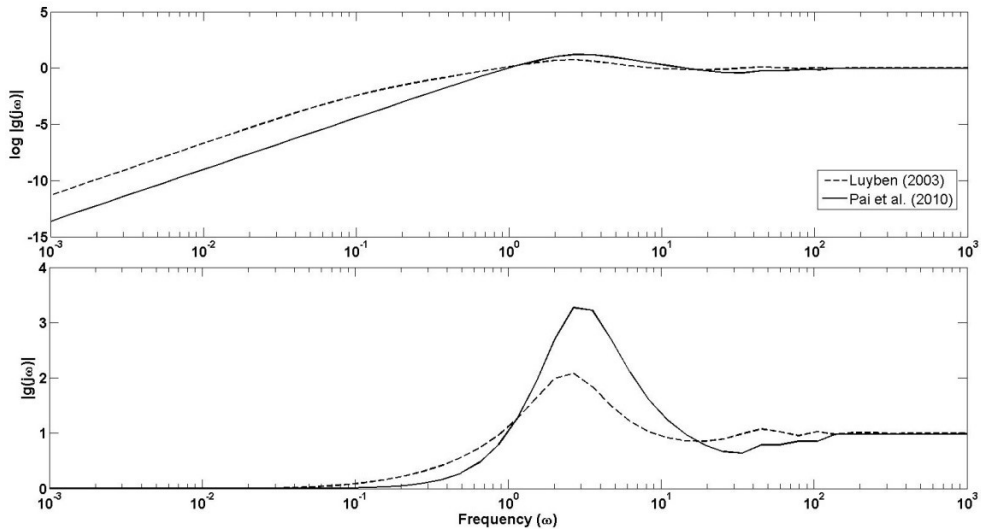

**Figure 2.** Graphical interpretations of the sensitivity function for example 1

The maximum absolute sensitivity $s_{\max}$ of the closed-loop with $G_{c1}(s)$ is less than that of the closed-loop with $G_{c2}(s)$. The term $s_{\max}$ is the measure of the largest amplification of disturbances. The relation $s_{\max} = \frac{1}{s_m}$ ($s_m$ is stability margin). For a closed-loop system, a better stability condition is known from the stability margin.

Thus, maximum absolute sensitivity, stability measure, and the robustness measure have a connection. The sensitivity performance indices in Table 1 indicate the superiority of the controller $G_{c1}(s)$ over $G_{c2}(s)$.

Table 1. Controller parameters and sensitivity performance indices for example 1

|  | Controller Parameters | | | Sensitivity indices (graphically calculated) | | | |
| --- | --- | --- | --- | --- | --- | --- | --- |
|  | $k_c$ | $\tau_I$ | $\tau_D$ | $s_{max}$ | $\omega_c$ | $\omega_{ms}$ | $\rho$ |
| Luyben (2003) | 1.69 | 11.5 | 1.15 | 2.09 | 0.8685 | 2.5 | 0.967 |
| Pai et al. (2010) | 4.06 | 2.68 | 0.65 | 3.33 | 0.9 | 3.089 | 0.794 |

**Example 2:** Autocatalytic Chemical Reactor

Consider an isothermal Continuous Stirred Tank Reactor (CSTR) effectuating an autocatalytic reaction $A + 2B \xrightarrow{k_1} 3B$ with a deactivation reaction $B \xrightarrow{k_2} C$ (Rao and Chidambaram, 2006), where the parameters $k_1$ and $k_2$ are reaction rates. A linearized model around the operating condition gives the transfer function model $G_p(s)$ with $k_p = -0.2679$, $\tau_0 = 41.6667$, $a_2 = 279.03$, $a_1 = -2.9781$ and $t_d = 10$, i.e.

$$G_p(s) = \frac{k_p(1-\tau_0 s)e^{-t_d s}}{a_2 s^2 + a_1 s + 1}, \quad a_2 > 0, a_1 < 0. \tag{20}$$

The presence of two unstable RHP *complex* poles in the plant transfer function $G_p(s)$ makes the problem complex. The controller structure $G_c(s)$ is a classic series PID with the additional lead-lag terminology. The controller parameters are $k_c = 1.3254$, $\tau_I = -86.251$, $\tau_D = 3.5807$, $\alpha_1 = 5$ and $\alpha_2 = 4.112$. Thus,

$$G_c(s) = k_c(1 + \frac{1}{\tau_I s} + \tau_D s)\frac{(\alpha_1 s + 1)}{(\alpha_2 s + 1)}. \tag{21}$$

With the above-mentioned parameters for the CSTR, the sensitivity function

$$g(j\omega) = \frac{\tau_I \omega^2(\alpha_2 a_2 \omega^2 - \alpha_2 - a_1) + j\tau_I \omega(1 - a_2 \omega^2 - \alpha_2 a_1 \omega^2)}{\tau_I \omega^2(\alpha_2 a_2 \omega^2 - \alpha_2 - a_1) + j\tau_I \omega(1 - a_2 \omega^2 - \alpha_2 a_1 \omega^2) + k_p k_c((\cos t_d \omega - j\sin t_d \omega)(a + jb))},$$

where

$$a = \tau_I \omega^2(-\alpha_1 - \tau_D + \tau_0 - \tau_0 \alpha_1 \tau_D \omega^2) + \tau_0 \alpha_1 \omega^2 + 1,$$
$$b = \tau_I \omega^3(-\alpha_1 \tau_D + \tau_0 \alpha_1 + \tau_0 \tau_D) + \omega(\alpha_1 + \tau_I - \tau_0).$$

Alternatively, in a brief and convenient form,

$$g(j\omega) = \frac{x_1 + jy_1}{k_p k_c(a\cos t_d \omega + b\sin t_d \omega) + x_1 + j(k_p k_c(b\cos t_d \omega - a\sin t_d \omega) + y_1)}, \tag{22}$$

where $x_1 = \tau_I \omega^2(\alpha_2 a_2 \omega^2 - \alpha_2 - a_1)$ and $y_1 = \tau_I \omega(1 - a_2 \omega^2 - \alpha_2 a_1 \omega^2)$. Consider an NMP zero associated with $G_p(s)$, i.e., $s = \frac{1}{\tau_0}$. The RHP open-loop poles of the plant transfer function $G_p(s)$ are

$$\alpha_1 = \frac{a_1 + \sqrt{a_1^2 - 4a_2}}{2a_2}, \quad \alpha_2 = \frac{a_1 - \sqrt{a_1^2 - 4a_2}}{2a_2}.$$

The pole-zero property of the CSTR sensitivity function has two non-minimum phase zeros and no RHP poles. Consider the singular point $s_0$ associated with the sensitivity integral, coincides with the NMP open-loop zero $s_0 = \frac{1}{\tau_0}$. The CSTR sensitivity integral, a consequence of Theorem 1 of the paper, becomes

$$\int_{-\infty}^{\infty} \log|g(j\omega)| \frac{\tau_0}{1 + \omega^2 \tau_0^2} d\omega = -\pi(\log\left|\frac{1 - \alpha_1 \tau_0}{1 + \overline{\alpha}_1 \tau_0}\right| + \log\left|\frac{1 - \alpha_2 \tau_0}{1 + \overline{\alpha}_2 \tau_0}\right|) \tag{23}$$

with $|g(j\omega)| = \dfrac{\sqrt{x_1^2 + y_1^2}}{\sqrt{x_1^2 + y_1^2 + k_p^2 k_c^2(a^2 + b^2) + 2k_p k_c((ax_1 + by_1)\cos t_d \omega + (bx_1 - ay_1)\sin t_d \omega)}}.$

By adopting appropriate variables for the CSTR specific case $g(j\omega)$ and using Theorem 2 of the paper, we get the CSTR logarithmic inequality in the Poisson-Jensen setting, i.e.

$$\log s_{\max} > \frac{\left|\frac{\pi}{2}\sum_{1\leq i\leq 2}\log\left|\frac{1-\alpha_i\tau_0}{1+\overline{\alpha}_i\tau_0}\right| + \log\max_{0\leq\omega\leq\omega_c}|g(j\omega)|\tan^{-1}\frac{\omega_c}{\sigma_k}\right|}{\frac{\pi}{2}-\tan^{-1}\frac{\omega_c}{\sigma_k}}, \quad \sum_{1\leq i\leq 2}\log\left|\frac{\sigma_k-\alpha_i}{\sigma_k+\overline{\alpha}_i}\right|<0. \tag{24a}$$

The above inequality (24a) is a consequence of equation (23) of the paper. Note that the parameters of equations (23) and (24a) have the same interpretations as the parameters of Theorem 1 of the paper. The CSTR sensitivity inequality can also be written in the Bode inequality setting as

$$\log s_{\max} \geq \underset{\omega_l \to \infty}{Lt} \frac{\left|-\pi\sum_{1\leq i\leq 2}\alpha_i + \omega_c\log\max_{0\leq\omega\leq\omega_c}|g(j\omega)|\right|}{\omega_l - \omega_c}, \tag{24b}$$

where $\sum_{1\leq i\leq 2}\alpha_i > 0$. The CSTR logarithmic inequality (24b) is a specific case of the Bode inequality (A.3).

It is observed from Figure 3 that, under a 10% model mismatch condition the sensitivity peak $s_{\max}$ goes on shifting upwards. This indicates more amplification of the disturbances beyond the frequency $\omega_c$. With the increase in the absolute maximum sensitivity $s_{\max}$, the stability margin becomes poorer indicating that the closed-loop system is losing its robustness, see Figure 3. Note that 10% uncertainties in gain, time constant and dead time are considered as a worse case of model mismatch scenarios.

Notably, the latter part of the $\log|g(j\omega)|$ trajectory in Figure 3 is oscillatory. The oscillatory behavior is attributed to the complex exponential term-coupled forward path transfer function with a larger plant dead time. The sensitivity performance indices, interpreted from the graph, are listed in Table 2. The logarithmic absolute maximum sensitivity of the CSTR is bounded above, i.e., $\log s_{\max} > 0.0766$ and $\log s_{\max} > 0.00335$ for the Poisson-Jensen and Bode settings respectively. The sensitivity lower bounds of the CSTR are calculated analytically in the Poisson-Jenson and Bode settings, see (24a) and (24b). The analytical and graphical calculations are listed in Table 2.

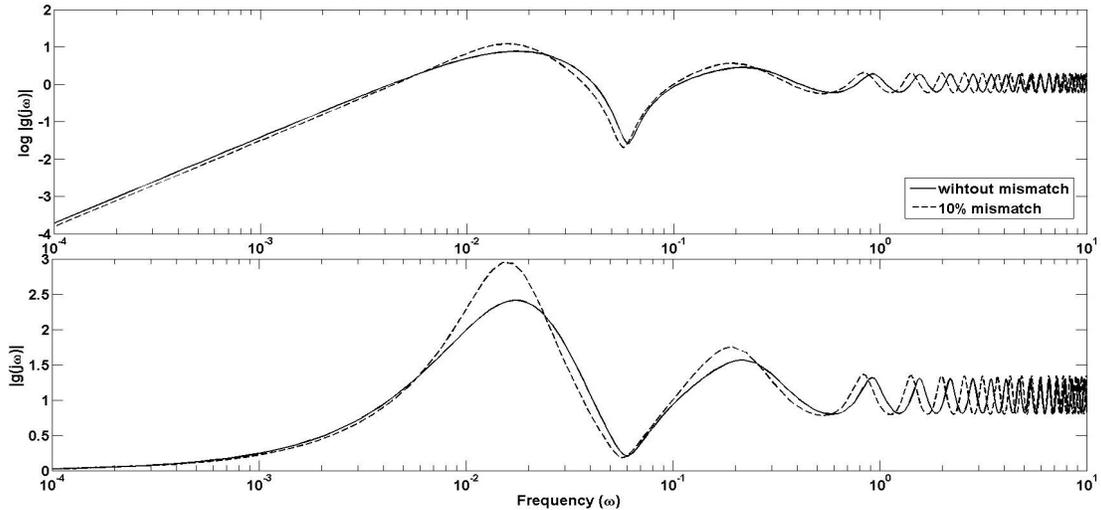

**Figure 3.** The sensitivity function for example 2 with 10% model mismatch

Table 2 reveals the following: (i) the logarithmic inequality in the Poisson-Jensen setting, accounts for the better lower bound in contrast to the Bode setting (ii) the graphical calculation of the absolute maximum sensitivity obeys the analytically calculated CSTR logarithmic inequality (iii) after the introduction of the 10% mismatch in the system parameters, the CSTR feedback system reveals a slightly better attenuation property for a given range of frequencies and slightly larger amplification.

Table 2. Sensitivity bounds and performance indices for example 2

| | Analytic sensitivity bound (log $s_{max}$) | | Sensitivity indices (graphically calculated) | | | |
|---|---|---|---|---|---|---|
| | Poisson-Jensen | Bode | $s_{max}$ | $\omega_c$ | $\omega_{ms}$ | $\rho$ |
| Example 2 without the model mismatch | >0.0766 | >0.00335 | 2.411 | 0.003972 | 0.0173 | 0.9333 |
| Example 2 with the 10% model mismatch | >0.0777 | >0.00336 | 2.947 | 0.004154 | 0.0152 | 0.9305 |

The feedback system retains the stability margin and preserves the robustness in the stability margin sense. The CSTR example reveals the usefulness of the logarithmic sensitivity integral and inequality Theorems of the paper.

**Example 3:** *SOPDT (Second Order Plus Dead Time)*

Consider a unity feedback system, where the transfer functions associated with the plant, controller, and sensitivity are $G_p(s)$, $G_c(s)$ and $g(s)$ respectively. The structures are

$$G_p(s) = \frac{k_p e^{-t_d s}}{(\tau_1 s - 1)(\tau_2 s + 1)}, \quad G_c(s) = k_c(1 + \frac{1}{\tau_I s} + \tau_D s), \quad g(s) = \frac{1}{1 + G(s)}, \quad (25)$$

where $G(s) = G_p(s) G_c(s)$. Furthermore,

$$g(s) = \frac{1}{1 + G_p(s) G_c(s)} = \frac{s\tau_I(\tau_1 s - 1)(\tau_2 s + 1)}{s\tau_I(\tau_1 s - 1)(\tau_2 s + 1) + k_p k_c (1 + \tau_I s + \tau_I \tau_D s^2) e^{-t_d s}},$$

$$g(j\omega) = \frac{\tau_I \omega^2 (\tau_2 - \tau_1) + j\tau_1 \omega(1 - \omega^2 \tau_1 \tau_2)}{\tau_I \omega^2 (\tau_2 - \tau_1) + j\tau_1 \omega(1 - \omega^2 \tau_1 \tau_2) + k_p k_c (1 + j\omega\tau_I - \tau_I \tau_D \omega^2) e^{-j t_d \omega}}$$

$$= \frac{x_2 + jy_2}{k_p k_c (x_3 \cos t_d \omega + \tau_I \omega \sin t_d \omega) + x_2 + j(k_p k_c (\tau_I \omega \cos t_d \omega - x_3 \sin t_d \omega) + y_2)}, \quad (26a)$$

where $x_2 = \tau_I \omega^2 (\tau_2 - \tau_1)$, $y_2 = \tau_1 \omega(1 - \omega^2 \tau_1 \tau_2)$, $x_3 = 1 - \tau_I \tau_D \omega^2$ and

$$|g(j\omega)| = \frac{\tau_I \omega \sqrt{\omega^4 \tau_1^2 \tau_2^2 + \omega^2(\tau_1^2 + \tau_2^2) + 1}}{\sqrt{x_2^2 + y_2^2 + k_p^2 k_c^2(x_3^2 + \tau_I^2 \omega^2) + 2k_p k_c ((x_3 x_2 + \tau_I y_2 \omega) \cos t_d \omega + (\tau_I x_2 \omega - x_3 y_2) \sin t_d \omega)}}. \quad (26b)$$

Consider a non-trivial *specific* plant transfer function $G_p(s)$ with an unstable RHP pole (Tan et al., 2002; Yang et al., 2002).

$$G_p(s) = \frac{e^{-0.939s}}{(5s - 1)(2.07s + 1)}. \quad (27)$$

equation (27) is a specific case, on the other hand, equation (25) is general. Numerous practical problems possess the form of an open-loop unstable system, e.g. chemical reactors with the exothermic reaction (Luyben, 1990), jacket-cooled CSTR (Luyben, 1998), etc. For this reason, equation (25) is examined via the introduction of sensitivity integrals into the case of equation (27). It is important to note that example 3 has an unstable open-loop pole. In the logarithmic integral, the singular point coincides with the NMP open-loop zero. The exponential delay $e^{-t_d s}$ annihilates the contribution of the closed path sensitivity integral in the Poisson-Jensen as well as the Bode settings, excluding the integral along the imaginary axis. Here, we consider three appealing controller structures available in the literature (Rao and Chidambaram, 2006; Shamsuzzoha and Lee, 2007, 2008) and examine them in a sensitivity inequality perspective. The first controller structure of Rao and Chidambaram (2006) is the same as that of two latter controller structures (Shamsuzzoha and Lee, 2007, 2008). The sensitivity inequality in Theorem 2 of the paper can be recast in the frameworks of the Poisson-Jensen and Bode sensitivity bounds. For the sensitivity function of equation (26b), the sensitivity integral, which is a consequence of Theorem 1 of the paper, becomes

$$\int_{-\infty}^{\infty} \log|g(j\omega)| \frac{\sigma}{\sigma^2 + \omega^2} d\omega = -\pi \log\left|\frac{\sigma - \alpha}{\sigma + \alpha}\right|, \quad (28)$$

where the terms $s_0 = \sigma \approx \dfrac{2}{t_d}$, $\alpha = \dfrac{1}{\tau_1}$ are associated with a system of equations (25) and (27). The singular point associated with the integral equation relates to the dead time of the plant. The closed-form solution to the integral equation is not possible, we combine equations (26a), (26b) and (28) with Theorem 2 of the paper. As a result of this, we get Poisson-Jensen sensitivity inequality, i.e.

$$\log s_{\max} > \dfrac{\left|\dfrac{\pi}{2}\log\left|\dfrac{\sigma - \alpha}{\sigma + \overline{\alpha}}\right| + \tan^{-1}\dfrac{\omega_c}{\sigma} \log\max_{0\leq\omega\leq\omega_c}|g(j\omega)|\right|}{\dfrac{\pi}{2} - \tan^{-1}\dfrac{\omega_c}{\sigma}}, \qquad (29a)$$

where $\log\left|\dfrac{\sigma - \alpha}{\sigma + \alpha}\right| < 0$. The logarithmic inequality (29a) can be regarded as an inequality version of equation (28). The sensitivity integral and the inequality in the Bode terminology becomes

$$\pi\alpha = \int_0^\infty \log|g(j\omega)|d\omega,$$

$$\log s_{\max} > \underset{\omega_l \to \infty}{Lt} \dfrac{\left|-\pi\alpha + \omega_c \log \max_{0\leq\omega\leq\omega_c}|g(j\omega)|\right|}{\omega_l - \omega_c}, \qquad (29b)$$

where the term $\alpha = \dfrac{1}{\tau_1}$ is associated with a system of equations (25) and (27). The Poisson-Jensen inequality gives a refined bound in contrast to the Bode, see Table 4. The controller structures of three appealing papers are the same in the general setting, see equation (25). However, they are different in the sense of different controller tuning parameters for the specific case, see Table 3. The difference is attributed to the different PID tuning methods. The same argument can be made about the difference between the three sensitivity integrals and three logarithmic inequalities that are associated with three controller structures (Rao and Chidambaram, 2006; Shamsuzzoha and Lee, 2007, 2008).

The controller efficacy is tested using the idea of sensitivity integrals and sensitivity inequalities for three well-established controller structures (Rao and Chidambaram, 2006; Shamsuzzoha and Lee, 2007, 2008). Figure 4 shows a plot of sensitivity functions for the three closed-loop systems. The controller of Shamsuzzoha and Lee (2008) has the maximum attenuation properties of disturbances in lower frequency ranges in contrast to two other controllers. However, the maximum value of $s_{\max}$, i.e. 1.608 occurring at the frequency $\omega_{ms} = 0.9655$, the controller of Shamsuzzoha and Lee (2008) will have the maximum amplification of disturbances occurring at the frequencies $\omega > \omega_c$.

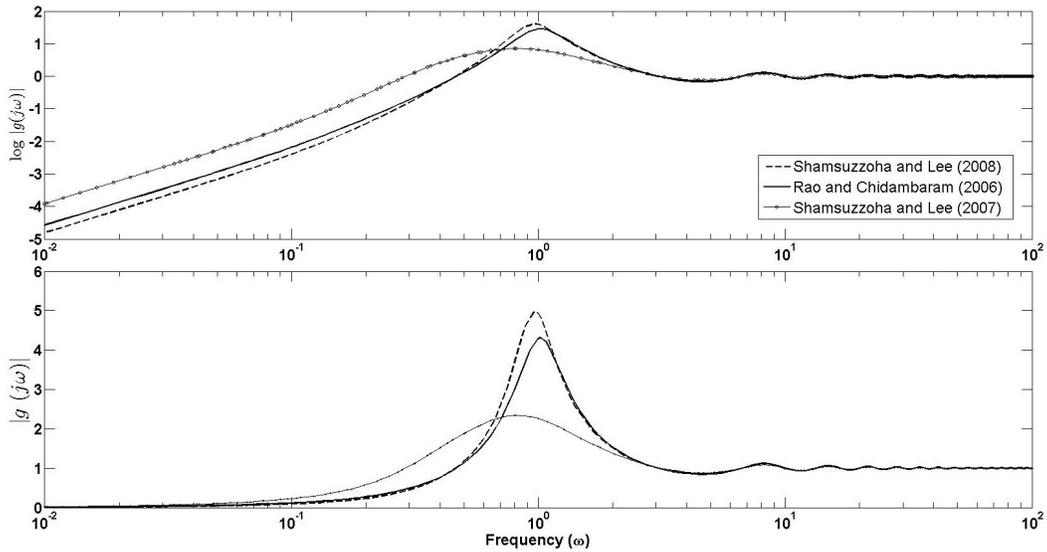

**Figure 4.** Graphical interpretations of the sensitivity function for example 3

Despite the minimum attenuation associated with Shamsuzzoha and Lee (2007) controller, it has the absolute maximum sensitivity $s_{max} = 2.338,$ which is relatively less in contrast to other controllers. That is indicative of a better stability margin, see Figure 4 and Table 3. The graphically calculated logarithmic absolute maximum sensitivity agrees with the analytically calculated sensitivity bounded above term, see Tables 3 and 4. Sensitivity bounds calculated analytically are listed in Table 4. That show agreements between the graphical results of the closed-loop systems and the sensitivity inequalities derived from Theorem 2 of the paper, see bounds (29a) and (29b).

**Table 3.** Controller parameters and sensitivity performance indices for example 3

| Method | Controller Parameters | | | Sensitivity indices, (graphically calculated) | | | |
|---|---|---|---|---|---|---|---|
| | $k_c$ | $\tau_I$ | $\tau_D$ | $s_{max}$ | $\omega_c$ | $\omega_{ms}$ | $\rho$ |
| Shamsuzzoha and Lee (2008) | 6.7051 | 5.4738 | 1.333 | 4.992 | 0.4498 | 0.9655 | 0.9443 |
| Rao and Chidambaram (2006) | 6.4285 | 6.4409 | 1.413 | 4.31 | 0.4479 | 1.023 | 0.9306 |
| Shamsuzzoha and Lee (2007) | 4.009 | 8.0327 | 1.6808 | 2.338 | 0.2805 | 0.8123 | 0.9536 |

Here, the loop robustness is adjudged by choosing the refined sensitivity lower bound and associating the bound with the logarithmic absolute maximum sensitivity in the inequality framework. The less logarithmic absolute maximum sensitivity reveals the higher stability margin. The higher margin preserves the robustness and brings apart from the fragility under the influence of input disturbances and plant parameter variations. Thus, the controller of Shamsuzzoha and Lee (2007) reveals greater robustness in contrast to other controllers.

**Table 4.** Sensitivity bounds for example 3

| Method | Analytic sensitivity bound ($\log s_{max}$) | |
|---|---|---|
| | Poisson-Jensen | Bode |
| Shamsuzzoha and Lee (2008) | >0.2267 | >0.0685 |
| Rao and Chidambaram (2006) | >0.2287 | >0.0691 |
| Shamsuzzoha and Lee (2007) | >0.2105 | >0.0641 |

*Modified sensitivity function for example 3*

To augment the ability of controllers, we recall Shamsuzzoha and Lee (2008) as an appealing example. The controller of Shamsuzzoha and Lee (2008) obeys equation (25) with specific plant parameters as well as controller parameters, see equation (27), and Table 3. To augment the ability of the controller subject to plant parameter variations, we associate two sensitivity functions, sensitivity function $g(s)$ and modified sensitivity function $\widetilde{g}(s)$. It is important to note that the sensitivity function $g(s)$ of Shamsuzzoha and Lee (2008) has an NMP zero and no RHP poles. Here in this paper, the modified sensitivity function $\widetilde{g}(s)$ will have no NMP zeros and no RHP poles, see equation (30a). The sensitivity integral and the inequality associated with Shamsuzzoha and Lee (2008) can be found in equations (28)-(29b) respectively. Here, we restrict our discussions to the modified sensitivity integral and related inequality associated with the forward path transfer function $\widetilde{G}(s)$ and the modified sensitivity function $\widetilde{g}(s)$. The modified sensitivity function

$$\widetilde{g}(s) = \frac{s\tau_I(\tau_1 s+1)(\tau_2 s+1)}{s\tau_I(\tau_1 s+1)(\tau_2 s+1) + k_p k_c (1+\tau_I s+\tau_I \tau_D s^2)e^{-t_d s}}, \quad (30a)$$

$$\widetilde{g}(j\omega) = \frac{-\tau_I \omega^2(\tau_2+\tau_1) - j\tau_I \omega(\tau_1 \tau_2 \omega^2-1)}{-\tau_I \omega^2(\tau_2+\tau_1) - j\tau_I \omega(\tau_1 \tau_2 \omega^2-1) + k_p k_c (1+j\omega\tau_I - \tau_I \tau_D \omega^2)e^{-j\omega t_d}}. \quad (30b)$$

Furthermore

$$\widetilde{g}(j\omega) = \frac{-x_4 - jy_4}{k_p k_c (x_5 \cos t_d \omega + \tau_I \omega \sin t_d \omega) - x_4 + j(k_p k_c (\tau_I \omega \cos t_d \omega - x_5 \sin t_d \omega) - y_4)}, \qquad (30c)$$

where $x_4 = \tau_I \omega^2 (\tau_2 + \tau_1)$, $y_4 = \tau_I \omega(\tau_1 \tau_2 \omega^2 - 1)$, $x_5 = 1 - \tau_I \tau_D \omega^2$ and

$$|\widetilde{g}(j\omega)| = \frac{\tau_I \omega^2 \sqrt{\tau_1^2 \tau_2^2 \omega^2 + \tau_1^2 + \tau_2^2}}{\sqrt{x_4^2 + y_4^2 + k_p^2 k_c^2 (x_5^2 + \tau^2 \omega^2) - 2 k_p k_c ((x_5 x_4 + y \tau \omega) \cos t_d \omega + (x \tau \omega - x_5 y) \sin t_d \omega)}}.$$

The modified sensitivity integral of Shamsuzzoha and Lee (2008) reduces to

$$\int_0^\infty \log|\widetilde{g}(j\omega)| \frac{\sigma}{\sigma^2 + \omega^2} d\omega = \frac{\pi}{2} \log|\widetilde{g}(\sigma_k)| = 0.$$

The above sensitivity integral is a specific case of Theorem 3 of the paper. The right-hand side of the above logarithmic sensitivity integral vanishes for the singular point $s_0 = \sigma \approx \frac{2}{t_d}$. The term $t_d$ is the dead time of the plant. The logarithmic inequality version of the above integral is more convenient to examine the controller for plant parameter variations. Thus,

$$\log s_{\max} > \frac{\left| \frac{\pi}{2} \log|\widetilde{g}(\sigma)| - \tan^{-1}\frac{\omega_c}{\sigma} \log \rho \right|}{\frac{\pi}{2} - \tan^{-1}\frac{\omega_c}{\sigma}} = \frac{\left| \tan^{-1}\frac{\omega_c}{\sigma} \log \rho \right|}{\frac{\pi}{2} - \tan^{-1}\frac{\omega_c}{\sigma}}.$$

The above inequality is a consequence of Theorem 4 of the paper. The term $\log \rho < 0$ has an interpretation of the attenuation properties of the input disturbances, and the term $\log s_{\max} > 0$ indicates the amplification of the input, see Theorem 4 of the paper. Figure 5 shows the sensitivity curves with a compensator and without compensator for example 3. Figure 5 unfolds the following: (i) the modified sensitivity function $\widetilde{g}(j\omega)$ has better attenuation properties in contrast to the sensitivity function $g(j\omega)$ for the frequency interval $0 \leq \omega \leq \omega_c$ (ii) the sensitivity peak associated with the modified sensitivity function $\widetilde{g}(j\omega)$ reduces, due to which the stability margin of the loop transfer function increases, which results in improved loop robustness. The modified sensitivity function $\widetilde{g}(j\omega)$ augments the robustness.

Furthermore, we study the controller of Shamsuzzoha and Lee (2008) by considering a set of three plant parameter variation cases. For the given three-parameter variation cases, the absolute maximum sensitivity, logarithmic sensitivity integral, and logarithmic inequalities are calculated. We consider the modified sensitivity function $\widetilde{g}(j\omega)$ accounting for the compensator. The same procedure is adopted by considering the sensitivity function $g(j\omega)$, see Figure 5 of the paper.

The model mismatching associated with the sensitivity function $g(j\omega)$ adds a more significant increase in the absolute maximum sensitivity which leads to a more significant decrease in the stability margin. Here, in the model mismatch case, 10% and 20% uncertainties in gain, dead time and time constant, are considered. The closed-loop system loses its robustness and becomes prone to disturbances. This may lead to system instability. On the other hand, the model mismatching associated with $\widetilde{g}(j\omega)$ adds a relatively less increase in absolute sensitivity. This leads to a relatively smaller decrease in the stability margin, see Figure 6.

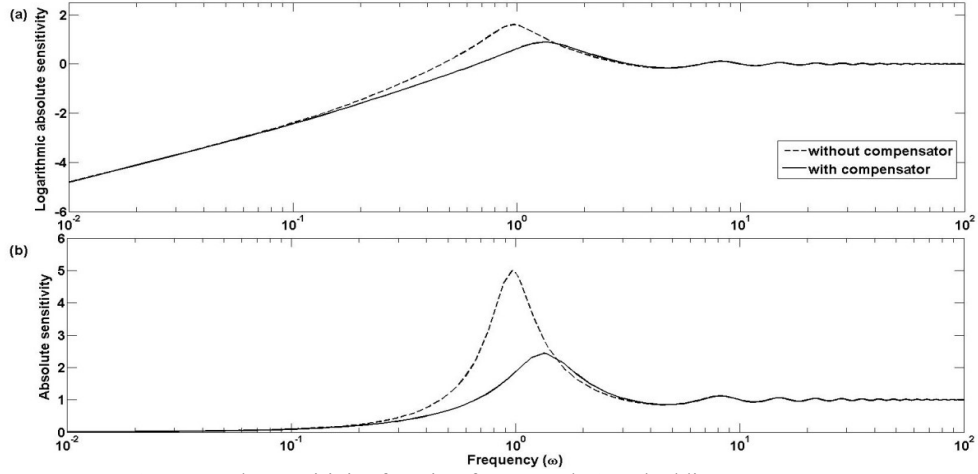

**Figure 5.** The sensitivity function for example 3 embedding compensator

For a 20 % mismatch, i.e. an increase in uncertainty, the compensated case displays a relatively less increase of 0.83 in the value of $s_{max}$. This ensures better stability margin. On the other hand, the uncompensated case displays a relatively higher increase in $s_{max}$, see Figure 6. Thus, there is a considerable improvement in the closed-loop performance due to the loop shaping. Table 5 presents graphically interpreted sensitivity performance indices with 10% and 20% uncertainties in the plant parameters.

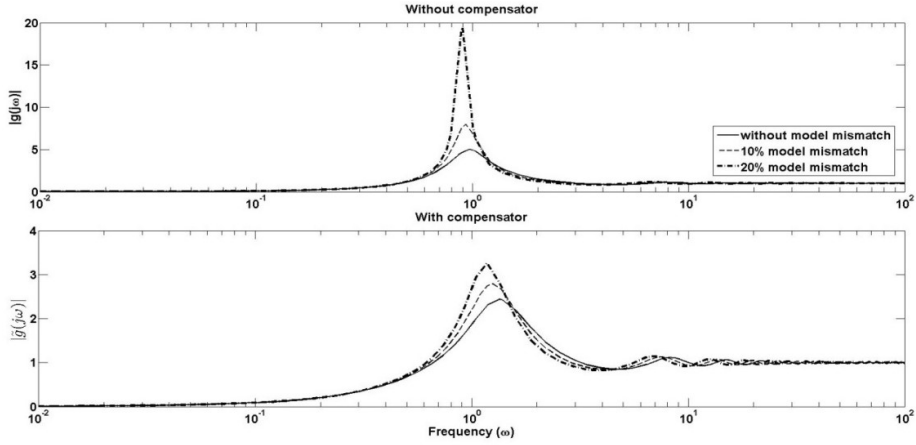

**Figure 6.** The sensitivity function for example 3 with model mismatching

Table 5 indicates that the modified sensitivity function $\tilde{g}(j\omega)$ is less prone to model mismatching effects, owing to less value of $s_{max}$. That suggests better robustness and better stability margin resulting from the loop shaping. On the other hand, the sensitivity function $g(j\omega)$ is relatively more prone to the model mismatching effects, which suggests less robustness and inadequate stability margin. Thus, the system may tend towards fragility. This demonstrates that Theorems 3 and 4 of the paper gives new insight into linear feedback systems in the robust vs. fragile sense.

**Table 5.** Sensitivity performance indices with compensator and uncertainties

| Uncertainty | Without compensator | | With compensator | |
|---|---|---|---|---|
| | $s_{max}$ | $\omega_{ms}$ | $s_{max}$ | $\omega_{ms}$ |
| Without mismatch | 4.992 | 0.9655 | 2.431 | 1.369 |
| 10% model mismatch | 7.972 | 0.9298 | 2.793 | 1.24 |
| 20% model mismatch | 19.45 | 0.8957 | 3.261 | 1.177 |

## 4. Conclusion

This paper sketches sensitivity integral Theorems and related logarithmic inequalities. The proofs of the Theorems combine the function's non-analyticity and complex integration. Furthermore, the analyticity notion of complex integration is utilized to augment the loop robustness. The sensitivity integrals and related lower bounds of the paper are formal robustness techniques. The modified sensitivity function theory of this paper is successfully demonstrated to augment the ability in the linear feedback system via the introduction of the compensator in the existing forward path. The results of the paper are phrased in the form of sensitivity integrals and sensitivity inequalities. This paper contributes towards 'filling the niche' between sensitivity integrals and chemical control systems by demonstrating the application of sensitivity integrals and inequalities, through Theorems of the paper, to the addressed chemical control systems. This paper also contributes towards the development of a formal theory of the sensitivity integral technique in 'detail' for the controller efficacy in contrast to the traditional, i.e. integral square error, integral absolute error, and total variation calculations. The technique of the paper is suggestive of application to other appealing control systems.

## Appendix

**Case 1:** Here, we rephrase the sensitivity integral of Theorem 1 in the Bode sensitivity integral equation setting. In the Bode integral setting, the singular point is not accounted for the closed path integral. On the other hand, the Poisson-Jensen integral equation accounts for the same. Since the Bode sensitivity integral is not a direct consequence of the Poisson-Jensen integral equation, we explain it succinctly. Recall the properties of the sensitivity functions $g(s)$ and $\tilde{g}(s)$ of Theorem 1 of the paper, we rephrase the sensitivity integral equation, i.e.

$$0 = \oint (\log g(s) + \sum_{1 \leq i \leq M} \log \frac{s+\alpha_i}{s-\alpha_i} + \sum_{1 \leq i \leq N} \log \frac{s-\beta_i}{s+\beta_i}) ds.$$

The above equation is a consequence of two results, i.e., the Green closed path integral and the analyticity. Furthermore, the closed path integral can be regarded as integral in the open RHP in the control framework. The closed path integral can be decomposed as

$$0 = j\int_{-\infty}^{\infty} \log g(j\omega) d\omega + \underset{R \to \infty}{Lt} \int_{\frac{\pi}{2}}^{-\frac{\pi}{2}} \log(g(Re^{j\gamma})) Re^{j\gamma} jd\gamma$$

$$+ \sum_{1 \leq i \leq M} j\int_{-\infty}^{\infty} \log \frac{j\omega+\alpha_i}{j\omega-\alpha_i} d\omega + \sum_{1 \leq i \leq M} \underset{R \to \infty}{Lt} \int_{\frac{\pi}{2}}^{-\frac{\pi}{2}} \log(\frac{Re^{j\gamma}+\alpha_i}{Re^{j\gamma}-\alpha_i}) j Re^{j\gamma} d\gamma$$

$$+ \sum_{1 \leq i \leq N} j\int_{-\infty}^{\infty} \log \frac{j\omega-\beta_i}{j\omega+\beta_i} d\omega + \sum_{1 \leq i \leq N} \underset{R \to \infty}{Lt} \int_{\frac{\pi}{2}}^{-\frac{\pi}{2}} \log(\frac{Re^{j\gamma}-\beta_i}{Re^{j\gamma}+\beta_i}) j Re^{j\gamma} d\gamma$$

$$= j\int_{-\infty}^{\infty} \log g(j\omega) d\omega + \underset{R \to \infty}{Lt} \int_{\frac{\pi}{2}}^{-\frac{\pi}{2}} \log(1+G(Re^{j\gamma})) Re^{j\gamma} jd\gamma$$

$$+ \sum_{1 \leq i \leq M} j\int_{-\infty}^{\infty} \log \frac{j\omega+\alpha_i}{j\omega-\alpha_i} d\omega + \sum_{1 \leq i \leq M} \underset{R \to \infty}{Lt} \int_{\frac{\pi}{2}}^{-\frac{\pi}{2}} \log(\frac{Re^{j\gamma}+\alpha_i}{Re^{j\gamma}-\alpha_i}) j Re^{j\gamma} d\gamma$$

$$+ \sum_{1 \leq i \leq N} j\int_{-\infty}^{\infty} \log \frac{j\omega-\beta_i}{j\omega+\beta_i} d\omega + \sum_{1 \leq i \leq N} \underset{R \to \infty}{Lt} \int_{\frac{\pi}{2}}^{-\frac{\pi}{2}} \log(\frac{Re^{j\gamma}-\beta_i}{Re^{j\gamma}+\beta_i}) j Re^{j\gamma} d\gamma. \qquad (A.1)$$

The poles of the sensitivity $g(s)$ function are attributed to the variable proportional gain $a$ of the feedback system. For a given bound on the variable gain, the RHP pole arises. The second term of the right-hand side accounts for the open-loop transfer function $G(s)$ with the variable gain $a$ and the unit relative degree, otherwise vanishes. Furthermore,

$$-a\pi = \underset{R\to\infty}{Lt} \int_{\frac{\pi}{2}}^{-\frac{\pi}{2}} \log(1+G(Re^{j\gamma}))Re^{j\gamma}\, jd\gamma,$$

$$\sum_{1\leq i\leq M}\int_{-\infty}^{\infty}\log\frac{j\omega+\alpha_i}{j\omega-\alpha_i}d\omega=0, \quad \sum_{1\leq i\leq N}\int_{-\infty}^{\infty}\log\frac{j\omega-\beta_i}{j\omega+\beta_i}d\omega=0,$$

$$-2\pi j\sum_i \alpha_i = \sum_{1\leq i\leq M}\underset{R\to\infty}{Lt}\int_{\frac{\pi}{2}}^{-\frac{\pi}{2}}\log(\frac{Re^{j\gamma}+\alpha_i}{Re^{j\gamma}-\alpha_i})j Re^{j\gamma}\, d\gamma,$$

$$2\pi j\sum_i \beta_i = \sum_{1\leq i\leq M}\underset{R\to\infty}{Lt}\int_{\frac{\pi}{2}}^{-\frac{\pi}{2}}\log(\frac{Re^{j\gamma}-\beta_i}{Re^{j\gamma}+\beta_i})j Re^{j\gamma}\, d\gamma.$$

After plugging the above set of relations into equation (A.1), we are led to

$$\frac{-a\pi}{2}+\pi\sum_{1\leq i\leq M}\alpha_i - \pi\sum_{1\leq i\leq N}\beta_i = \int_0^{\infty}\log|g(j\omega)|d\omega. \tag{A.2}$$

Equation (6) and (A.2) are two different sensitivity integral equations. Their usefulness, i.e., the Poisson-Jensen and Bode, can be studied by consideration of appealing practical problems. For the NMP open-loop zero case, associated with the open-loop transfer function $G(s)$, the sensitivity integral equation assumes an appealing structure, see Equation (7). This setup is a convenient form to compute sensitivity inequality.

**Case 2:** The significant difference between the Poisson-Jensen and Bode terminologies is attributed to the singular point lying in the right half-plane as well as associated with the complex integration. The Poisson-Jensen accounts for the singularity in the integral, on the other hand, the Bode does not account for the singularity. Suppose the linear feedback system does not have the real open-loop NMP zeros, we consider the arbitrary singularity lying in the right-half place in the Poisson-Jensen setting, on the other hand, the singularity does not hold for the Bode setting.

Consider a general form of the Bode integral

$$\int_0^{\infty}\log|g(j\omega)|d\omega = \frac{-a\pi}{2}+\pi\sum_{1\leq i\leq M}\alpha_i - \pi\sum_{1\leq i\leq N}\beta_i.$$

Its inequality can be stated as

$$\omega_c \log\rho + (\omega_l-\omega_c)\log s_{max} \geq \frac{-a\pi}{2}+\pi\sum_{1\leq i\leq M}\alpha_i - \pi\sum_{1\leq i\leq N}\beta_i.$$

Furthermore

$$(\omega_l-\omega_c)\log s_{max} \geq \frac{-a\pi}{2}+\pi\sum_{1\leq i\leq M}\alpha_i - \pi\sum_{1\leq i\leq N}\beta_i - \omega_c\log\rho.$$

After invoking the condition $\frac{-a\pi}{2}+\pi\sum_{1\leq i\leq M}\alpha_i-\pi\sum_{1\leq i\leq N}\beta_i > 0$,

$$\log s_{max} \geq \underset{\omega_l\to\infty}{Lt}\frac{\left|\frac{a\pi}{2}-\pi\sum_{1\leq i\leq M}\alpha_i+\pi\sum_{1\leq i\leq N}\beta_i+\omega_c\log\rho\right|}{\omega_l-\omega_c}. \tag{A.3}$$

This inequality is general and covers a wide class of linear feedback systems. This holds, irrespective of the presence of the open-loop NMP zeros. For the relative degree greater than one or presence of the delay term, the term $\frac{a\pi}{2} = 0$. For the stable linear feedback system, the term $\beta_i$ vanishes. Note that instability in the linear feedback system is attributed to a range of values of the open-loop gain $a$.

**Case 3:** Here, we rephrase the sensitivity integrals of Theorem 3 of the paper in the Bode terminology, i.e.

$$\int_0^\infty \log|\tilde{g}(j\omega)|d\omega = \frac{-a\pi}{2}.$$

Note that the variable gain $a$ is associated with the modified forward path transfer function $\tilde{G}(s)$ in lieu of $G(s)$. The modified sensitivity function in the Bode terminology and the useful functions are

$$\tilde{g}(s) = \kappa(s)g(s), \; \kappa(s) = \prod_{1 \leq i \leq M}(\frac{s+\alpha_i}{s-\alpha_i}) \prod_{1 \leq i \leq N}(\frac{s-\beta_i}{s+\beta_i}),$$

$$\tilde{G}(s) = \frac{1-\kappa(s)+G(s)}{\kappa(s)}.$$

For the unit relative degree of $\tilde{G}(s)$, the above integral has a non-vanishing property, otherwise vanishes.

**Case 4:** Consider a general form of the Bode integral for the sensitivity function $\tilde{g}(s)$,

$$\int_0^\infty \log|\tilde{g}(j\omega)|d\omega = \frac{-a\pi}{2},$$

where the gain $a$ is associated with the modified forward path transfer function $\tilde{G}(s)$. Its inequality can be stated as

$$\omega_c \log \rho + (\omega_l - \omega_c) \log s_{\max} \geq \frac{-a\pi}{2}.$$

Furthermore

$$(\omega_l - \omega_c) \log s_{\max} \geq \frac{-a\pi}{2} - \omega_c \log \rho.$$

By invoking the condition $\frac{-a\pi}{2} - \omega_c \log \rho > 0$,

$$\log s_{\max} \geq \underset{\omega_l \to \infty}{Lt} \frac{|\frac{a\pi}{2} + \omega_c \log \rho|}{\omega_l - \omega_c}.$$

For the relative degree greater than one or the presence of the delay term, the term $\frac{a\pi}{2} = 0$ and the inequality simplify to $\log s_{\max} \geq \underset{\omega_l \to \infty}{Lt} \frac{|\omega_c \log \rho|}{\omega_l - \omega_c}$. Note that the notations $\rho$ and $s_{\max}$ of the case 4 are associated with Theorem 4 of the paper.